\documentclass{aims}
\usepackage{amsmath}
\usepackage{amssymb}
  \usepackage{paralist}
 \usepackage[colorlinks=true]{hyperref}
\hypersetup{urlcolor=blue, citecolor=red}

\def\currentvolume{X}
 \def\currentissue{X}
  \def\currentyear{200X}
   \def\currentmonth{XX}
    \def\ppages{X--XX}

\newtheorem{theorem}{Theorem}[section]

\theoremstyle{definition}

\title[Estimates for Fourier integrals]
      {Dispersive type estimates for Fourier integrals and applications to hyperbolic systems}

\author[Michael Ruzhansky and Jens Wirth]{}


 \email{m.ruzhansky@imperial.ac.uk}
 \email{j.wirth@imperial.ac.uk}

\thanks{The authors are supported by EPSRC grant  EP/E062873/1.}

\begin{document}
\maketitle

% Enter the first author's name and address:
\centerline{\scshape Michael Ruzhansky and Jens Wirth}
\medskip
{\footnotesize
% please put the address of the first author
 \centerline{Department of Mathematics, Imperial College London}
   \centerline{180 Queen's Gate}
   \centerline{London, SW7 2AZ, UK}
} % Do not forget to end the {\footnotesize by the sign }

\bigskip

% The name of the associate editor will be entered by an editorial staff
 \centerline{(Communicated by XXX)}

%The abstract of your paper
\begin{abstract}
In this note we provide dispersive estimates for Fourier integrals with parameter-dependent phase functions in terms of geometric quantities of associated families of Fresnel surfaces. The results are based on a multi-dimensional van der Corput lemma due to the first author. 

Applications to dispersive estimates for hyperbolic systems and scalar higher order hyperbolic equations are also discussed.
\end{abstract}

%The title of your section 1
\section{Introduction}
Dispersive estimates for solutions to linear evolution equations are vitally important for the study of a wide range of problems, including questions of stability of solutions to corresponding non-linear equations, descriptions of large-time asymptotic profiles of solutions or scattering results based on the time-dependent approach.

For the particular situation of hyperbolic evolution equations with a hyperbolic time-asymptotic profile, such estimates can be based on corresponding decay estimates for Fourier integral operators with real phases and certain well-controled amplitude functions. The aim of this note is to provide some robust estimates for such Fourier integrals in terms of geometric properties of associated Fresnel surfaces. 

%The title of your section 2
\section{Fourier integral operators, notation and main results}
We are interested in $L^1$--$L^\infty$ estimates for a family of Fourier integral operators depending on a parameter $t\in\mathbb R_+$ (thought of as time later on) 
\begin{equation}\label{eq:Tdef}
   T_t u (x) = \int_{\mathbb R^n} \mathrm e^{\mathrm i ( x\cdot\xi + t\phi(t,x,\xi) )}a(t,x,\xi)\, \widehat u(\xi) \mathrm d\xi,
\end{equation}
where $a(t,x,\xi)$ is uniformly in $t$ a (global) symbol of order zero,
\begin{equation}\label{eq:symbEstXi}
  |\partial_\xi^\alpha a(t,x,\xi) | \le C_\alpha \langle\xi\rangle^{-|\alpha|},
\end{equation}
with $\langle\xi\rangle=\sqrt{1+|\xi|^2}$, and $\phi(t,x,\xi)$ is a homogeneous real-valued function satisfying 
\begin{equation}
  C^{-1} |\xi| \le \phi(t,x,\xi) \le C |\xi|,\qquad |\partial_\xi^\alpha \phi(t,x,\xi)| \le C_\alpha' |\xi|^{1-|\alpha|}
\end{equation}
uniform in $t$ and $x$. For simplicity, we always assume $u\in\mathcal S(\mathbb R^n)$, estimates extend by continuity arguments to larger spaces. If we are only interested in estimates of $T_tu$ in the $L^\infty$-norm, we can treat $t$ and $x$ both as parameters and derive estimates for the corresponding stationary model operator
\begin{equation}\label{eq:Model}
  T : u \mapsto \int_{\mathbb R^n} \mathrm e^{\mathrm i(x\cdot\xi + \phi(\xi)) } a(\xi) \widehat u(\xi) \mathrm d\xi.
\end{equation}
 A precise understanding of the appearing constants in the estimate translates to estimates of the $t$-behaviour of the full model \eqref{eq:Tdef}. 

\subsection{Fresnel surfaces and contact indices} We turn to the analysis of \eqref{eq:Model} and define the geometric quantities which determine its $L^\infty$ behaviour. The consideration follows essentially Sugimoto \cite{sugi94}, \cite{Sug96}, for the non-degenerate situations see the earlier works of Strichartz \cite{Str70} or Brenner \cite{Bre75}. 

We assume that the phase $\phi(\xi)$ is non-negative, 1-homogeneous and smooth, and define the associated Fresnel surface (also called slowness surface)
\begin{equation}
  \Sigma = \{ \xi\in\mathbb R^n : \phi(\xi) = 1 \}.
\end{equation}
We introduce two indices 
for such a Fresnel surface $\Sigma$, assuming that
it is of class $C^k$ with $k$ being sufficiently large. 
For $p\in\Sigma$ we denote by $\mathbb T_p$ the  
tangent hyperplane to $\Sigma$ at $p$. Then for any plane $H$ 
of dimension 2 which contains $p$ and the normal of
$\Sigma$ at $p$ we denote by $\gamma(\Sigma,p,H)$ 
the order of contact between the curve
$\Sigma\cap H$ and its tangent $H\cap \mathbb T_p$. Furthermore, we set
\begin{equation}
 \gamma(\Sigma) = \sup_{p} \sup_{H} \gamma(\Sigma,p,H),\qquad 
 \gamma_0 (\Sigma) = \sup_{p} \inf_H \gamma(\Sigma, p,H).
\end{equation}
The definition implies directly that
$2\le\gamma_0(\Sigma)\le\gamma(\Sigma)$. 
Furthermore, for isotropic problems $\Sigma$ is a dilation of the sphere $\mathbb S^{n-1}$ and 
$\gamma(\mathbb S^{n-1})=\gamma_0(\mathbb S^{n-1})=2$.
Moreover, if the Gaussian curvature of $\Sigma$ never
vanishes $\Sigma$ is convex and we have $\gamma(\Sigma)=2$.

These indices are directly related to mapping properties of \eqref{eq:Model}. The following theorem is a combination of results from \cite{sugi94}, \cite{Sug96} in the analytic setting,
improved by the author \cite{Ru09}, \cite{RS10} to the smooth 
setting as well as to the limited regularity cases.
Moreover, the corresponding results in 
\cite{Ru09}, \cite{RS10} do not require the positive homogeneity
of the phase in $\xi$ and allow the dependence on parameters.

As usual, $B^r_{p,q}(\mathbb R^n)$ denotes the Besov space of regularity $r$ modelled over $L^p$.
\begin{theorem}
  Assume that $\phi$ is smooth on $\mathbb R^{1+n}\times(\mathbb R^n\setminus\{0\})$, positive and $1$-homogeneous in $\xi$, $\Sigma$ smooth  and $a\in S^0_{1,0}(\mathbb R^n)$. Then the operator $T$ defined by 
  \eqref{eq:Model} maps 
  \begin{equation}
  B^r_{1,2}(\mathbb R^n) \to L^\infty(\mathbb R^n)
  \end{equation}
  where $r = n-\frac{n-1}{\gamma(\Sigma)}$ if $\Sigma$ is convex and 
  $r = n-\frac{1}{\gamma_0(\Sigma)}$ if $\Sigma$ is arbitrary. 
\end{theorem}
The operator norm in this estimate depends on finitely many symbol seminorms of $a$, the estimates of $\phi$ and its derivatives together with some quantitative measures of the contact between $\Sigma$ and tangent lines. 

To control the order of contact of $\Sigma$ by tangent lines quantitatively, we will introduce another quantity $\varkappa(\Sigma)$. First, assume that $\Sigma$ is convex and that it is of class $C^{\gamma(\Sigma)+1}$. For $p\in\Sigma$, rotating $\Sigma$ if
necessary, we may assume that it is parameterised
by points $\{(y,h(y)), y\in \Omega\}$ near $p$ for an open set
$\Omega\subset\mathbb R^{n-1}$. For $p=(y,h(y))$, let us define
\begin{equation}
\varkappa(\Sigma,p)=\inf_{|\omega|=1} \sum_{j=2}^{\gamma(\Sigma)}
\left| \frac{\partial^j}{\partial\rho^j} h(y+\rho\omega)|_{\rho=0}
\right|.
\end{equation}
From the definition of $\gamma(\Sigma)$
it follows that $\varkappa(\Sigma,p)>0$ for all $p\in\Omega$.
Indeed, from the definition of $\gamma(\Sigma,p,H)$ it follows
that if $\omega$ is such that $y+\rho\omega\in H$, then
\begin{equation}\label{EQ:ass-kappa1}
\varkappa(\Sigma,p,H)=\left|\frac{\partial^{\gamma(\Sigma,p,H)}}
{\partial\rho^{\gamma(\Sigma,p,H)}} h(y+\rho\omega)|_{\rho=0}
\right|>0.
\end{equation} 
Now, we clearly have
$\sum_{j=2}^{\gamma(\Sigma)}
\left| \frac{\partial^j}{\partial\rho^j} h(y+\rho\omega)|_{\rho=0}
\right|\geq \varkappa(\Sigma,p,H)$,
and hence we have $\varkappa(\Sigma,p)>0$ since the set
$|\omega|=1$ is compact.
Noticing that $\varkappa(\Sigma,p)$
is a continuous function of $p$, by compactness of
$\Sigma$ it
follows that if we define
\begin{equation}
\varkappa(\Sigma)=\min_{p\in\Sigma} \varkappa(\Sigma,p),
\end{equation}
then $\varkappa(\Sigma)>0$.
If $\Sigma$ is not convex, we define similarly
\begin{equation}
\varkappa_0(\Sigma)=\min_{p\in\Sigma} 
\sup_{|\omega|=1} \sum_{j=2}^{\gamma_0(\Sigma)}
\left| \frac{\partial^j}{\partial\rho^j} h(y+\rho\omega)|_{\rho=0}
\right|.
\end{equation}
Again, we have $\varkappa_0(\Sigma)>0$.

Considering again the example of spheres $\Sigma=r\mathbb S^{n-1}$ of radius $r$, we see that $\gamma(\Sigma)=\gamma_0(\Sigma)=2$ and that $\varkappa(\Sigma)=\varkappa_0(\Sigma)\sim \frac 1r$ is a measure of scalar curvature. 

\subsection{Main results}\label{sec:main}
Estimates of \eqref{eq:Tdef} depend on indices of a family of Fresnel surfaces defined in terms 
$\phi(t,x,\xi)$,
\begin{equation}\label{eq:Sigma-Def}
   \Sigma_{t,x}  =\{ \xi\in\mathbb R^n : \phi(t,x,\xi) = 1\},\qquad t\ge t_0, x\in\mathbb R^n.
\end{equation}
For a first result we assume that $\Sigma_{t,x}$ is convex for all $t\ge t_0$ and all $x\in\mathbb R^n$.

\begin{theorem}\label{thm:2.3}
Assume that $\phi(t,x,\xi)$ is real-valued, continuous in $t$ and $x$ and smooth in $\xi\in\mathbb R^n\setminus\{0\}$, homogeneous of order one in $\xi$ and such that for some $t_0>0$ and $C>0$ we have
\begin{equation}
    C^{-1} |\xi| \le \phi(t,x,\xi) \le C|\xi| \quad\text{and}\quad |\partial_\xi^\alpha \phi(t,x,\xi)|\le C|\xi|^{1-|\alpha|}
\end{equation}
for all $t\ge t_0$, all $x$, $\xi\ne0$ and all multi-indices $\alpha$ with $|\alpha|\le \max\{ \gamma+1, \lfloor (n-1)/\gamma\rfloor +2\}$. Assume further that all the sets $\Sigma_{t,x}$ defined by \eqref{eq:Sigma-Def}
are convex for all $t\ge t_0$ and assume that 
\begin{equation}
 \limsup_{t\to\infty}\sup_x \gamma(\Sigma_{t,x}) \le\gamma
 \quad\text{together with}\quad
 \liminf_{t\to\infty} \inf_x \varkappa(\Sigma_{t,x})>0.
\end{equation}
Suppose further that the amplitude $a(t,x,\xi)$ satisfies 
\begin{equation}
   |\partial_\xi^\alpha a(t,x,\xi)| \le C_\alpha \langle\xi\rangle^{-|\alpha|}\quad\text{for all}\quad |\alpha|\le \lfloor (n-1)/\gamma\rfloor +1.
\end{equation}
Then the operator family $T_t$ defined by \eqref{eq:Tdef} satisfies for all $t\ge t_0$ the estimate
\begin{equation}
   \| T_t u\|_{L^\infty(\mathbb R^n)} \le C t^{-\frac{n-1}\gamma} \|u\|_{B^r_{1,2} (\mathbb R^n)}
\end{equation}
with $r=n-\frac{n-1}\gamma$.
\end{theorem}

In the non-convex situation decay rates can be much lower. Results can be improved by using a more detailed analysis of the geometric situation near points of high contact order.

\begin{theorem}\label{thm:2.4}
Assume that $\phi(t,x,\xi)$ is real-valued, continuous in $t$ and $x$ and smooth in $\xi\in\mathbb R^n\setminus\{0\}$, homogeneous of order one in $\xi$ and such that for some $t_0>0$ and $C>0$ we have
\begin{equation}
    C^{-1} |\xi| \le \phi(t,x,\xi) \le C|\xi| \quad\text{and}\quad |\partial_\xi^\alpha \phi(t,x,\xi)|\le C|\xi|^{1-|\alpha|}
\end{equation}
for all $t\ge t_0$, all $x$, $\xi\ne0$ and all multi-indices $\alpha$ with $|\alpha|\le\gamma_0+1$. Assume further that all the sets $\Sigma_{t,x}$ defined by \eqref{eq:Sigma-Def} satisfy 
\begin{equation}
 \limsup_{t\to\infty}\sup_x \gamma_0(\Sigma_{t,x}) \le\gamma_0
 \quad\text{together with}\quad
 \liminf_{t\to\infty} \inf_x \varkappa_0(\Sigma_{t,x})>0.
\end{equation}
 Suppose further that the amplitude $a(t,x,\xi)$ satisfies 
\begin{equation}
   |\partial_\xi^\alpha a(t,x,\xi)| \le C_\alpha 
   \langle\xi\rangle^{-|\alpha|}\quad\text{for all}\quad 
   |\alpha|\le 1.
\end{equation}
Then the operator family $T_t$ defined by 
\eqref{eq:Tdef} satisfies for all $t\ge t_0$ the estimate
\begin{equation}
   \| T_t u\|_{L^\infty(\mathbb R^n)} \le C 
   t^{-\frac{1}{\gamma_0}} \|u\|_{B^r_{1,2} (\mathbb R^n)}
\end{equation}
with $r=n-\frac{1}{\gamma_0}$.
\end{theorem}

The proofs of both theorems are similar to the ones in \cite{RW10} and full details are omitted here. The dependence on $x$ is parametric, and the assumption assure that the resulting estimates are uniform in $x$.

If we make further regularity assumptions with respect to the $x$-variable it follows that the operators $T_t$ are uniformly $L^2$-bounded such that interpolation implies the usual $L^p$--$L^q$ decay estimates. For completeness, we include a corresponding statement. The assumptions are chosen such that after substituting $x=ty$ and 
$t\xi=\eta$ Theorem~2.5 of \cite{RS06} applies:

\begin{theorem}[\cite{RS06}]\label{Th2.5}
Let $\Gamma_y,\Gamma_\xi\subset\mathbb R^n$ be open cones.
Let operator $T$ be defined by
\begin{equation}\label{T2}
 Tu(x)
 =\int_{\mathbb R^n}\int_{\mathbb R^n} 
 e^{i(x\cdot\xi+\varphi(y,\xi))}a(y,\xi)u(y)\ \mathrm dy 
 \ \mathrm d\xi,
\end{equation}
where $a(y,\xi)\in C^\infty(\mathbb R^n_y\times\mathbb R^n_\xi)$,
${\rm supp}\ a\subset \Gamma_y\times\Gamma_\xi$, and
$\varphi(y,\xi)\in C^\infty(\Gamma_y\times\Gamma_\xi)$ is
a real-valued function.
Assume that
\[
|\partial_y^\alpha\partial_\xi^\beta a(y,\xi)|
\leq C_{\alpha\beta},
\]
for $|\alpha|, |\beta|\leq 2n+1$.
Also assume that
\[
|\det\partial_y\partial_\xi\varphi(y,\xi)|\geq C>0
\quad\text{on}\quad\Gamma_y\times\Gamma_\xi,
\]
that
$$
|\partial_\xi\varphi(x,\xi)-\partial_\xi
\varphi(y,\xi)|\geq C|x-y|\quad
\text{for} \quad
x,y\in\Gamma_y,\; \xi\in\Gamma_\xi,
$$
and that
\[
|\partial_y^\alpha\partial_\xi \varphi(y,\xi)\leq C_{\alpha}|,
\qquad
|\partial_y\partial_\xi^\beta \varphi(y,\xi)|\leq C_{\beta}
\quad\text{on supp}\; a
\]
for $1\leq |\alpha|,|\beta|\leq 2n+2$.
Then the operator $T$ is $L^2(\mathbb R^n)$-bounded, and satisfies
\[
\|T\|_{L^2\to L^2}\leq C
\;\text{\rm sup}_{|\alpha|,|\beta|\leq 2n+1}
 \|\partial_y^\alpha\partial_\xi^\beta a(y,\xi)\|
  _{L^\infty(\mathbb R^n_y\times\mathbb R^n_\xi)}.
\]
\end{theorem}

We can now formulate the corresponding $L^2$-result:

\begin{theorem}\label{THM:L2}
Let 
\begin{equation}\label{eq:Tdef11}
   T_t u (x) = \int_{\mathbb R^n} \mathrm e^{\mathrm i ( x\cdot\xi + t\phi(t,x,\xi) )}a(t,x,\xi)\, \widehat u(\xi) \mathrm d\xi,
\end{equation}
Assume $\phi(t,x,\xi)$ is real-valued, $1$-homogeneous in $\xi$ and satisfies
\begin{equation}\label{eq:phase-asss}
   |\det(\mathrm I + t \partial_x\partial_\xi \phi(t,x,\xi))| \ge C_0 >0,\qquad |t^{|\alpha|}\partial_x^\alpha\partial_\xi^\beta \phi(t,x,\xi)| \le C_{\alpha,\beta}
\end{equation}
uniform in $t\ge t_0$, $x$, $\xi\ne0$ and for $1\le |\alpha|,|\beta|\le  2n+2$. Assume further that $a(t,x,\xi)$ is supported in 
$t |\xi|\ge C$ for some constant $C$ and that
\begin{equation}
   |t^{|\alpha|}\partial_x^\alpha\partial_\xi^\beta a(t,x,\xi)| \le C_{\alpha,\beta}% t^{-|\alpha|} %|\xi|^{-|\beta|}
\end{equation}
for all $|\alpha|, |\beta|\le 2n+2$ and $t\ge t_0$. Then the family $T_t$ in \eqref{eq:Tdef11}, $t\ge t_0$, is uniformly $L^2$-bounded,
\begin{equation}
   \|T_t u\|_{L^2(\mathbb R^n)} \le C \| u\|_{L^2(\mathbb R^n)}.
\end{equation} 
\end{theorem}
Let us sketch how Theorem \ref{Th2.5} implies 
Theorem \ref{THM:L2}. Using homogeneity, we can rewrite 
\eqref{eq:Tdef11} as
\begin{equation}\label{eq:Tdef2}
   T_t u (x) = t^{-n}\int_{\mathbb R^n} \mathrm e^{\mathrm i ( 
   z\cdot\eta + \phi_t(t,z,\eta) )}a_t(t,z,\eta)\, 
   \widehat u(\eta/t) 
\mathrm d\eta,
\end{equation}
with $x=tz, t\xi=\eta$, and
$$
\phi_t(t,z,\eta)=\phi(t,tz,\eta),\quad
a_t(t,z,\eta)=a(t,tz,\eta/t).
$$
By rescaling, the powers of $t$ in the $L^2$-estimate
cancel, and it can be shown that the adjoint to the operator in
\eqref{eq:Tdef2} satisfies assumptions on Theorem
\ref{Th2.5} uniformly for $t\geq t_0$.

We also note that the assumptions on the mixed derivatives of 
$\phi$ as in \eqref{eq:phase-asss} can be relaxed by looking only at
$|\alpha|=1$ or $\beta=1$, as in Theorem \ref{Th2.5}.

\subsection{Appendix: Fourier transform of surface carried measures} Let $\Sigma\subseteq\mathbb R^n$ be a closed hypersurface. Then the indices $\gamma(\Sigma)$ and $\gamma_0(\Sigma)$ are related to decay estimates for Fourier transforms of surface carried measures. Consider for $f\in C^\infty(\Sigma)$ the integral
\begin{equation}
  u(x) = \int_\Sigma \mathrm e^{\mathrm i x\cdot\xi} f(\xi) \mathrm d\xi,
\end{equation}
where $\mathrm d\xi$ denotes the $(n-1)$-dimensional surface measure on $\Sigma$. Then the following theorem is valid.
\begin{theorem}
\begin{enumerate}
\item Assume $\Sigma$ is a smooth hypersurface and convex. Then
\begin{equation}
   |u(x)| \le C \langle x\rangle^{-\frac{n-1}{\gamma(\Sigma)}} \|f\|_{C^r(\Sigma)}
\end{equation}
with $r\ge (n-1)/\gamma(\Sigma)+1$.
\item For general (non-convex) smooth hypersurface $\Sigma$ it follows that
\begin{equation}
   |u(x)| \le C \langle x\rangle^{-\frac{1}{\gamma_0(\Sigma)}} \|f\|_{C^1(\Sigma)}.
\end{equation}
\end{enumerate}
\end{theorem}
The proof of the first statement in \cite{sugi94} uses real analyticity of the surface, the stronger results are due to methods of \cite{Ru09} and \cite{RS10}. Smoothness of the surface can be further relaxed to the class $C^{\gamma}$ with $\gamma=\gamma(\Sigma)$ (or $\gamma_0(\Sigma)$, respectively).

\section{Applications to hyperbolic evolution equations}
The estimates discussed so far have applications for dispersive estimates of hyperbolic evolution equations. Without aiming for completeness we will discuss some possible applications.% and show some improvements of older results.

\subsection{Wave models} As a first example we consider wave models with time-dependent coefficients,
\begin{equation}\label{eq:wave}
   \partial_t^2 u - \sum_{i,j=1}^n a_{ij}(t) \partial_{x_i}\partial_{x_j} u% + 2b(t) \partial_t u + c(t) u
    = 0
\end{equation}
in $\mathbb R\times\mathbb R^n$ and with Cauchy data $u(0,\cdot)=u_0$, $u_t(0,\cdot)=u_1$, both assumed to be Schwartz for simplicity. Denoting
\begin{equation}
  \mathcal T\{\ell\} = \{ a\in C^\infty(\mathbb R) : |\partial_t^k a(t) | \le c_k \langle t\rangle^{\ell-k} \},
\end{equation}
we assume  $a_{ij} \in \mathcal T\{0\}$, $a_{ij}(t)=\overline{a_{ji}(t)}$, together with the strict hyperbolicity assumption
\begin{equation}
a(t,\xi) =  \sum_{i,j=1}^n a_{ij}(t) \xi_i\xi_j \ge C |\xi|^2,
\end{equation}
for some $C>0$ uniform in $t$. Within the zone $\{(1+t)|\xi|\gtrsim 1\}$ solutions to \eqref{eq:wave} are given by Fourier integrals of the form \eqref{eq:Tdef} with phase functions
\begin{equation}
  \phi_\pm(t,\xi) =\pm  \frac 1t \int_0^t \sqrt{a(\theta,\xi)}\mathrm d\theta.
\end{equation}
In this situation the Fresnel surfaces are given by
\begin{equation}
  \Sigma_t = \left\{ \xi\in\mathbb R^n :  t = \int_0^t \sqrt{a(\theta,\xi)}\mathrm d\theta \right\}
\end{equation}
and satisfy $\gamma(\Sigma_t)=2$ together with uniform upper and lower bounds for $\varkappa(\Sigma_t)$. 
This follows directly from \cite[Lemma~3.1]{R04}. Thus Theorem~\ref{thm:2.3} is applicable and
in combination with a uniform bound on the $L^2$-energy one concludes that solutions to \eqref{eq:wave} satisfy
\begin{equation}
   \|\partial_t u(t,\cdot)\|_q + \| \nabla u(t,\cdot)\|_q \le C_{pq} (1+t)^{-\frac{n-1}2\left(\frac1p-\frac1q\right)} \left( 
   \|u_1\|_{W^{p,r_p+1}} + \|u_2\|_{W^{p,r_p}}\right)
\end{equation}
with $1<p\le 2\le q<\infty$, $pq=p+q$ and $r_p = n(1/p-1/q)$.

For a detailed consideration of such models we refer to Reissig \cite{R04}, \cite{RS05a}. It is possible to include further lower order terms $b(t)\partial_t u$, $\sum_j c_j(t) \partial_{x_j} u$ and $d(t)u$ with $b, c_j\in \mathcal T\{-1\}$ and $d\in\mathcal T\{-2\}$ (together with further assumptions to control the low-frequency behaviour). We refer to \cite{Wir06} for an example of such a treatment. If lower order terms belong to different classes, the asymptotic behaviour of solutions might change its type, see, e.g., \cite{Wir07} a discussion about dissipative wave equations and effective dissipation.

The assumptions on coefficients can be weakened to classes
\begin{equation}
  \mathcal T_\nu\{\ell\} = \{ a\in C^\infty(\mathbb R) : |\partial_t^k a(t) | \le c_k \langle t \log (e+t) ^\nu\rangle^{\ell-k} \}
\end{equation}
with parameter $\nu\in[0,1]$ where $\nu\in(0,1)$ implies a small loss of decay and $\nu=1$ yields a polynomial loss. It can be shown that without further assumptions $\nu>1$ leads generically to supra-polynomial growth of the energy. Additional assumptions to compensate this behaviour for wave equations with variable speed of propagation were considered by Hirosawa \cite{Hiro07}, see also \cite{Wir10b}, but the estimates of Theorems~\ref{thm:2.3} and \ref{thm:2.4} need to be refined in order to obtain dispersive estimates in this situation.

\subsection{Higher order scalar equations}\label{sec:3.2}
In a similar way we can treat homogeneous higher order equations of the form
\begin{align}
   &\mathrm D_t^m u + \sum_{k=0}^{m-1} \sum_{|\alpha|=m-k} a_{k,\alpha} (t) \mathrm D_t^k \mathrm D_x^\alpha u = 0,\\
   & \mathrm D_t^j u(0,\cdot) = u_j,\quad j=0,1,\ldots, m-1,
\end{align}
where $\mathrm D = -\mathrm i\partial$, $a_{k,\alpha}\in\mathcal T\{0\}$ and the characteristic equation
\begin{equation}
  \tau^m + \sum_{k=0}^{m-1}\sum_{|\alpha|=m-k} a_{k,\alpha} (t) \tau^k \xi^\alpha = 0
\end{equation}
has $m$ uniformly distinct real solutions $\tau_1(t,\xi)$,\ldots, $\tau_m(t,\xi)$, $|\tau_i(t,\xi) - \tau_j(t,\xi) | \ge C|\xi|$ for $i\ne j$ and some $C>0$. Under the stronger assumption that $a_{k,\alpha}'(t)\in L^1(\mathbb R_+)$, this model was treated by Matsuyama and the first author, \cite{MR10}, the full form follows from \cite{RW10}. If we assume
\begin{equation}\label{EQ:hh-int}
 %  \max_{1\le j\le m} \sup_{T\ge 0}\sup_{\xi\ne0}
    \bigg| \sum_{k\ne j} \int_0^t \frac{\partial_t \tau_j(\theta,\xi)}{\tau_j(\theta,\xi)-\tau_k(\theta,\xi)}\mathrm d \theta  \bigg|<C
\end{equation}
uniform in $j$, $t$ and $\xi\ne0$, solutions are given as sums of Fourier integrals of the form \eqref{eq:Tdef} with phases
\begin{equation}\label{eq:phi-j}
   \phi_j(t,\xi) = \frac1t \int_0^t \tau_j(\theta,\xi)\mathrm d\theta
\end{equation}
and uniformly bounded amplitudes. Again, the results from Section~\ref{sec:main} provide dispersive estimates for solutions to such problems. We remark that for higher order equations the associated Fresnel surfaces $\Sigma_t$ need not be convex and therefore the full generality of Theorems~\ref{thm:2.3} and \ref{thm:2.4} is of importance here.

We note that the analysis in \cite{MR10} is based on the
asymptotic integration method for constructing solutions
for equations with homogeneous symbols,
and it requires only $C^1$-regularity of coefficients,
i.e. only that $a_{k,\alpha}\in C^1$ in terms of the regularity
in time, which is crucial in applications to the Kirchhoff
equations, see e.g. \cite{MR09}.

\subsection{Hyperbolic systems}\label{sec:3.3} 
In \cite{RW08}, the authors announced 
results on dispersive estimates 
for a general class of pseudo-differentiable 
hyperbolic systems with time-dependent coefficients 
and arbitrary $\gamma$.
D'Abbico, Lucente and Taglialatela, \cite{dAbbico:2009}, 
considered differential 
hyperbolic systems with $\gamma=2$ , while in \cite{RW10} 
the authors gave a rigorous
treatment for general pseudo-differentiable hyperbolic 
systems with time-dependent coefficients and 
arbitrary $\gamma$. For the sake of simplicity we 
formulate the results for the homogeneous situation. We consider
\begin{equation}\label{eq:system}
   \mathrm D_t U = A(t,\mathrm D_x) U,\qquad U(0,\cdot)=U_0
\end{equation}
for a $m\times m$ matrix $A(t,\mathrm D_x)$ of $t$-dependent Fourier multipliers. We denote
\begin{equation}
  \mathcal S\{m_1,m_2\} = \{ a(t,\xi) \in C^\infty(\mathbb R\times (\mathbb R^n\setminus\{0\}) ) : 
   |\xi|^{|\alpha|-m_1} \mathrm D_\xi^\alpha a(t,\xi) \in \mathcal T\{m_1\} \},
\end{equation}
the condition being uniform in $\xi$.
We assume that the symbol $A(t,\xi)\in\mathcal S\{1,0\}$ is $1$-homogeneous in $\xi$, a self-adjoint matrix and its eigenvalues
$\tau_1(t,\xi)$,\ldots, $\tau_m(t,\xi)$ are uniformly distinct, $|\tau_i(t,\xi) - \tau_j(t,\xi) | \ge C|\xi|$ for $i\ne j$ and some $C>0$. Without further assumptions it follows that the solutions to \eqref{eq:system} can be represented in the form
\eqref{eq:Tdef} with phases given by \eqref{eq:phi-j}.

Assumptions can be relaxed, without assuming self-adjointness we have to require real eigenvalues together with a uniform symmetrisability of the matrix $A(t,\xi)$ together with further assumptions to control the large-time influence of lower order terms. We refer to \cite{RW10} for details.

\medskip
% The data information below will be filled by AIMS editorial staff
Received xxxx 20xx; revised xxxx 20xx.
\medskip

\end{document}